%Paper: q-alg/9501022
%From: Charilaos Aneziris <aneziris@ifh.de>
%Date: Wed, 18 Jan 1995 06:13:14 +0100 (MET)
%Date (revised): Thu, 19 Jan 1995 03:35:23 +0100 (MET)
%Date (revised): Fri, 20 Jan 1995 03:40:16 +0100 (MET)
%Date (revised): Fri, 20 Jan 1995 05:14:23 +0100 (MET)

\magnification=\magstep1
\parskip 0pt
\baselineskip 16pt
\null
\vskip 1 cm
\centerline {\bf IS A KNOT CLASSIFICATION POSSIBLE?}
\bigskip
\centerline {Charilaos Aneziris}
\centerline {\it Institut f\"ur Hochenergienphysik Zeuthen}
\centerline {\it DESY Deutsches Elektronen-Synchrotron}
\centerline {\it Platanenallee 6}
\centerline {\it 15738 Zeuthen}
\bigskip
\centerline {\bf Abstract}
\bigskip
The goal of this paper is to discuss the possibility of finding an algorithm
that can give all distinct knots
up to a desired complexity. Two such algorithms are presented, one based on
projections on a plane, the other on closed self-avoiding walks.
\bigskip
\centerline {1. \qquad INTRODUCTION}
\bigskip
Knot theory has been one of the most thoroughly studied and researched areas of
Mathematical Physics during the last and the current decade. While work on
knots
has been carried out almost exclusively by mathematicians for about a century,
during recent years a number of unexpected applications in science, from
topological field theories in High Energy Physics [1] to DNA properties in
Biology [2]
has given an enormous boost to the subject and has expanded the number of
disciplines related to it. In particular, it is known that Wilson lines are
the observables in topological field theories. As Witten has shown in Ref. 1,
when the gauge group is $SU(2)$,
the correlation functions of Wilson lines obey relations similar
to the ones satisfied for the Jones polynomials of the corresponding knots and
links, which were initially introduced in knot
theory. A knot classification would therefore be extremely useful in studying
the observables of topological field theories and their correlation functions.
\par \medskip
In spite, however, of the work done on knot theory and the recent progress that
has been achieved, we have yet to come up with some method or formula that if
appropriately applied would yield a complete table with all the distinct types
of knot. Current knot tables that are in use have either been found by hand or
through incomplete algorithms that fail when the knot complexity becomes too
large; it is hardly surprising, therefore, that such tables often differ with
each
other. One might think that the task of obtaining such tables should have been
facilitated by the discovery of a number of knot and link characteristics
ranging from the Alexander polynomials [3] to the Vassiliev invariants
[4]. While such characteristics are indeed very helpful in distinguishing
inequivalent knots, they fail to provide the full
picture. This is not only due to the
fact that it is possible for inequivalent knots to have the same polynomial,
but
in addition, given a polynomial, one cannot always tell which knot, if any,
corresponds to this polynomial.
\par \medskip
Until just a few years ago the best knot tables existing were based on the ones
provided by Reidemeister and Seifert in the 1930's [5], which did not
proceed
further than the first 250 simplest knots. A few errors creeped in these tables
however; it was only in 1974, for example, that Perko proved that the knots
$10_{161}$ and $10_{162}$ of Rolfsen's list are the same [6]. In 1985
Thistlethwaite
came up with a table of about 13,000 knots [7],but the possibility of
duplications
remains and the problem is still far from being resolved. Similarly, a number
of attempts to prove that the program is intrinsically solvable were either
incomplete [8], or impractical [9].
\par \medskip
In this paper we first provide the necessary definitions and background
and review some of the relevant work carried out in the past. Then we develop
two
possible algorithms for knot classification. The first one is based on
two-dimensional projections and bears some similarity to
Thistlethwaite's work, while the other on the theory of self-avoiding walks.
\bigskip
\centerline {2. \qquad DEFINITIONS}
\bigskip
A {\bf knot} is defined as a one-to-one homeomorphism from $S^1$ to some
manifold $\cal M$.
\par \medskip
According to the definition above, knots may not possess ``double" points; the
figure ``8" for example is not a knot. Such figures however are sometimes
called
``singular knots" and play an important role in defining the Vassiliev
invariants of (non-singular) knots.
\par \medskip
In addition, the definition above is too ``broad", and it obviously allows for
a continuously infinite number of distinct knots; no algorithm may ever hope to
classify the set of these knots. On the other side, however, if one considers
as equivalent all knots isotopic to each other, one gets a set that is too
``narrow", since all knots are isotopic to $S^1$ and therefore to each other.
This is because intrinsically all knots are the same; only when they are
embedded in some target space one may possible look for distinct classes of
knots. Knot equivalence is therefore defined as follows.
\par \medskip
Two knots $K_1$ and $K_2$ are called {\bf ambient isotopic} and considered
equivalent iff there is a homeomorphism $f$ from $[0,1]$ to the set of knots
such that $f(0)=K_1$ and $f(1)=K_2$.
\par \medskip
This definition may also be stated as follows. Let $g_1$ and $g_2$ be the
functions from $S^1$ to $\cal M$ defining $K_1$ and $K_2$. The two knots are
ambient isotopic iff a function $f(t,s)$ from $[0,1] \times S^1$ to $\cal M$
exists which is continuous with respect to both $t$ and $s$ and satisfies
$f(0,s)=g_1(s)$, $f(1,s)=g_2(s)$, and
$\forall t \in [0,1]: s_1 \ne s_2 \Rightarrow f(t,s_1) \ne f(t,s_2)$.
\par \medskip
Proving that ambient isotopy is an equivalence relation is trivial; one may
easily check that $K \sim K$, $K_1 \sim K_2 \Rightarrow K_2 \sim K_1$ and
$K_1 \sim K_2 \wedge K_2 \sim K_3 \Rightarrow K_1 \sim K_3$, where $\sim$
denotes ambient isotopy. Therefore ambient isotopy divides knots into distinct
classes; it is the classification and enumeration of such classes that we are
going to study from now on.
\par \medskip
Two subtle points shall be mentioned
now. First, if one reverses
orientation, that is replace $f(s)$ with $f(-s)$ as the defining function, one
gets a different knot, even though these two knots consist of exactly the same
points. The same thing happens of course by an orientation preserving
reparametrization. When orientation is preserved, however, the two knots are
always ambient
isotopic, while when orientation is reversed, they need not be so.
Similarly, when one replaces $f(s)$ with $-f(s)$ one gets the {\it mirror
image} of $f(s)$; a knot may or may not be ambient isotopic to its mirror
image.
If it is, it is called {\bf amphicheiral}, while if it is not, it is called
{\bf chiral}. In most of the paper, however, knots will be considered
equivalent to both their mirror images and to inversely oriented knots, while
in section 11 the problem of how to find which knots are invertible and
amphicheiral will be raised once more.
\par \medskip
Finally, before leaving this section, a few examples for
various target spaces shall be provided, while from section 3 the discussion
will be entirely directed towards solving
the ${\cal M}=S^3$ case.
\par \medskip
For ${\cal M}=R^1$ obviously no knots are allowed, while for ${\cal M}=S^1$ or
${\cal M}=R^2$ there are just two distinct classes of knots, one that is
clockwise and one that is counterclockwise, which are related through
orientation
inversion. For ${\cal M}=S^2$ these two classes coincide. On the other side,
for ${\cal M}=S^n$ and ${\cal M}=R^n$ for $n \ge 4$, all knots are equivalent
to each other and just one class exists.
\par \medskip
A more interesting case arises for the torus, ${\cal M}=T^2=S^1 \times S^1$.
While all homotopic knots are ambient isotopic and vice versa, not all homotopy
classes are
suitable to generate knots. Indeed a class $(m,n)$ may only generate knots if
the Greatest Common Denominator (GCD) of $m$ and $n$ is $1$. If it is larger
than 1,
all loops belonging to such class will possess ``double" points and thus do
not fit the definition of
knots. The set of knot classes is therefore $\{(m,n):m,n \in Z \wedge {\rm GCD}
(m,n)=1\}$.
\par \medskip
The ${\cal M}=S^3$ case, however, is even more interesting since it still
remains unsolved, and it is this case to be discussed from now on.
\bigskip
\centerline {3. \qquad KNOTS ON $S^3$; THE REIDEMEISTER MOVES}
\bigskip
While for ${\cal M} = T^2$ distinct knot classes exist due to the
non-triviality of the fundamental group, $S^3$ possesses a trivial fundamental
group
and thus all loops are homotopic; non-trivial knots do exist, however, because
although they can be continuously deformed to identity, such a process involves
self-intersecting loops. It is impossible, for example, to deform the
{\it trefoil} (Figure 1a) to the {\it unknot} (Figure 1b) without passing
through a self-intersecting loop.
\par \medskip
The question of knots in $S^3$ was drastically simplified and advanced during
the 1920's when Reidemeister proved that knots may be studied through their
projections in $S^2$ [10]. Given a certain knot, one may replace it with a
two-dimensional projection such that no more than two knot points may be
projected
on the same point of $S^2$, while the neighborhood of all double points should
look
as in Figure 2a (and not as in Figure 2b). The second condition can be more
rigorously stated as follows. Let $f(s_1)$ and $f(s_2)$ be two points of $S^3$
that are projected at the same point of $S^2$. The difference of the
derivatives
at these two points, $f'(s_2)-f'(s_1)$, should not be orthogonal to
$S^2$. In addition, one must keep track of
overcrossings and undercrossings so that Figure 3a is considered different from
3b. Projections satisfying these conditions are called {\bf regular}.
\par \medskip
Once regular projections of two knots are given, the two knots are equivalent
if one may continuously deform one projection to the other. Even however if
such a continuous deformation is not possible, the two knots are equivalent if
the projections are connected through a sequence of moves such as the ones
shown in Figures 4a, 4b and 4c. These three kinds of moves are called the
{\bf Reidemeister} moves [11]. It is important to notice that no double
points other
than the ones shown in Figure 4 are allowed in the bold face rectangulars, and
subsequently the move of Figure 5 is not allowed. In addition, the part of the
knot outside the bold face rectangular should not change during the move.
\par \medskip
While one may easily notice that by performing any of these three moves, or any
combination of them, one gets a knot equivalent to the initial one, the
important result obtained by Reidemeister is the converse, meaning that these
three kinds of moves are sufficient; if two knot projections cannot be
connected through a sequence of just these three moves, the knots are
definitely inequivalent. This result will be crucial when devising an
algorithm to get distinct classes of inequivalent knots.
\bigskip
\centerline {4. \qquad NAMING A KNOT}
\bigskip
A knot, as defined earlier, is a one-to-one homeomorphism from $S^1$ to $S^3$.
No algorithm can of course give an actual knot, since knots are continuous
objects while algorithms are discrete. It is therefore essential to provide a
knot class with a discrete ``name", that is a finite grouping of integer
numbers which uniquely characterizes a knot class; a desirable algorithm is
then going to be a computer program whose input corresponds to some desired
knot complexity, while the output is a list of knot class ``names" whose
complexity does not exceed the input's value. (Such an input is necesary
because the number of knot classes is
infinite and without this input the program would run forever. The definition
of a knot complexity will depend of course on the algorithm, and therefore
two such definitions shall be given, one for each algorithm to be
discussed, while an interesting point which shall be very briefly touched
here is
whether a knot that is more complex to another with respect to one algorithm,
may be less complex with respect to another).
\par \medskip
Let us now consider a two-dimensional projection of a knot such as the one
shown in Figure 6 satisfying the rules set in section 3. First one chooses a
starting
point on the projection and one of the two orientations; in Figure 6 this point
is O and the orientation is clockwise at the neighborhood of O. As one
travels starting from O along the chosen orientation on the knot
projection, one assigns the double points encountered, successive natural
numbers
starting from $1$ and going upwards, until one returns to O. Each double point
is thus acquiring two such numbers, one for the overcrossing segment and one
for the undercrossing. Thus for example the double points of Figure 6 are
assigned the following numbers.
\vskip 12 cm \centerline {T A B L E \quad I} \medskip
\vbox{\catcode`\*=\active \def*{\hphantom{0}}
\offinterlineskip
\halign{\strut#&\vrule#\quad&\hfil#\hfil&\quad
\vrule#\quad&\hfil#\hfil&\quad\vrule#\quad&\hfil#\hfil&\quad\vrule#\cr
\noalign{\hrule}
&&\omit\bf \quad Double Point \quad && \bf Overcrossing Number &&
\bf Undercrossing Number
&\cr \noalign {\hrule}
&&A&& *1&& 12& \cr &&B&& *2&& *7& \cr &&C&& 10&& *3& \cr &&D&& *4&& *9& \cr
&&E&& *5&& 14& \cr &&F&& *6&& 15& \cr &&G&& *8&& 11& \cr &&H&& 16&& 13& \cr
\noalign{\hrule}}}
\par \medskip Each double point is now assigned a pair whose left
and right elements are the overcrossing and the undercrossing numbers
respectively. The name of the projection is the set of all such pairs. For
the case of Figure 6, for example, the name of the knot is $\{ (1,12),(2,7),
(10,3),(4,9),(5,14),(6,15),(8,11),(16,13) \}$. The reason for such a choice is
that while permuting the order of pairs does not affect the projection,
permuting the order of the two numbers of a pair does.
\par \medskip Once a knot is assigned a name using the process above, or any
other process for that reason, there are four questions one needs to answer
before going ahead.
\par \medskip 1) Can all knot classes be assigned such a name?
\par \smallskip 2) Is such a name unique for a certain class?
\par \smallskip 3) Is there always a knot class to be assigned for any
possible name?
\par \smallskip 4) If so, is such a class unique?
\par \medskip The answers to these questions come out to be: yes, no, no, yes.
As a result of the negative answers to the second and third questions, two
problems arise.
\par \medskip 1) How can one proceed to choose just one name for a particular
class?
\par \smallskip 2) How can one distinguish the names that do correspond to
knot classes from the ones that do not?
\par \medskip In the next two sections these two
questions are going to be discussed.
First, however, two points should be mentioned.
One is the notion of knot complexity; here it is going to be the
number of double points.
Therefore once some natural number
$n$ is set as input, one expects to get as output all knot classes,
named the way described in this section, that have at least one projection
with not more than $n$ double points.
\par \medskip The other is about the number of possible names.
While for any
natural number $n$ this number is $(2n)! \over n!$ since
there are $(2n)!$ ways to order the numbers from $1$ to $2n$ and $n!$ ways to
order the $n$ pairs, the actual number of knot classes is much smaller. Before
proceeding it would be interesting to give a look at the following table.
\bigskip
\centerline {T A B L E \quad II} \medskip
\vbox{\catcode`\*=\active \def*{\hphantom{0}}
\offinterlineskip
\halign{\strut#&\vrule#\quad&\hfil#\hfil&\quad
\vrule#\quad&\hfil#\hfil&\quad\vrule#\quad&\hfil#\hfil&\quad\vrule#\cr
\noalign{\hrule}
&&\omit \qquad $n$ \qquad && \qquad $(2n)! \quad \div \quad n!$ \qquad &&
\bf Actual \quad Number \quad of \quad Knots &\cr \noalign
{\hrule}
&&0&& **********1&& *1&\cr &&1&& **********2&& *0&\cr
&&2&& *********12&& *0&\cr
&&3&& ********120&& *1&\cr
&&4&& *******1680&& *1&\cr &&5&& ******30240&& *2&\cr &&6&& *****665280&&
*3&\cr
&&7&& ***17297280&& *7&\cr
&&8&& **518918400&& 21&\cr &&9&& 17643225600&& 49&\cr
\noalign{\hrule}}}
\bigskip
\centerline {5. \qquad CHOOSING A NAME}
\bigskip
While in section 4 a procedure for naming knot projections was given, by doing
so one is still far from reaching the ultimate goal. In order to get distinct
knot
classes one has to make sure that no knot class appears more than once under
different names. To achieve now this goal one should first know which names
correspond to equivalent knots and why, and then provide a criterion that
will eliminate all names but one, so that the surviving name will appear in the
output as the sole representative of its knot class.
\par \medskip There are two reasons for multiple names to be assigned to one
and the same
knot class. First, by varying the starting point O of section 4 and/or
inverting orientation one gets an identical projection but possibly under a
different name. Second, one may additionally perform Reidemeister moves on a
certain projection and obtain different projections, which belong however
to the same knot class.
In each of these two cases a name will be affected as follows.
\par \medskip By varying the starting point while preserving orientation, all
pairs $(i,j)$ become $(k+i,k+j)$ for some $k \in \{ 0,1,2,...,2n-1 \}$. If
however orientation is reversed, all pairs $(i,j)$ become instead $(k-i,k-j)$.
In both cases the value of $k$ depends on the location of the new starting
point. One may  notice that the numbers appearing in the new pairs do
not necessarily belong to $\{ 1,2,...,2n \}$. In such a case they have to be
replaced by mod$2n$. Such a case is very likely to occur in all name changes
to be studied from now on, and the substitution of $r$ by $r \rm {mod}
2n$ shall always be implied without further saying.
\par \medskip One may therefore establish the following rule. Two names $N_1$
and $N_2$, which of course contain the same number of pairs, belong to the
same projection if and only if there is a number $k \in \{ 0,1,2,
...,2n-1 \}$ such that either $(i,j) \in N_1 \Rightarrow (k+i,k+j) \in N_2$,
or $(i,j) \in N_1 \Rightarrow (k-i,k-j) \in N_2$. Having now established the
criterion above, all one has to do now is find the names of some projection
which are at most $4n$, give some rule that will eliminate all names but one,
and keep this one name as the sole representative of the projection. This
name will later be compared to names of other projections of the same knot
class so that eventually only one name should survive and appear at the output.
It is important to notice here that this is only possible because the maximum
possible number of equivalent names, $4n$, is {\it finite}.
\par \medskip In order now for some rule that compares names to be valid, one
must first make sure that for any two names different from each other this
rule keeps exactly one of them, if a name happens to be compared to
itself it must survive, and finally if $N_1$ is eliminated by $N_2$ and $N_2$
by $N_3$, then $N_1$ has to be eliminated by $N_3$. While there is a number of
criteria that can satisfy these demands, the one to be chosen, the
``lexicographical", shall be described at the end of section 6, after the
discussion of the problem of ``impossible" names.
\par \medskip One last point will be mentioned here before moving on to
the Reidemeister moves. Let $N_1$ and $N_2$ be two names containing the same
number of pairs such that $(i,j) \in N_1 \Rightarrow (j,i) \in N_2$. Such
names may either belong to the same knot class (Figure 7a) or to mirror
images (Figure 7b). According to what was said in section 2, in either case
$N_1$ and $N_2$ will be considered identical for the time being. In order to
avoid duplication we shall only keep names where $1$ appears at the left side
of a pair. If after some procedure we get a name where $1$ appears at the
right, then all pairs $(i,j)$ shall be replaced by $(j,i)$ so that $1$ appears
at the left.
\par \medskip When one performs a ``first" Reidemeister move (Figure 4a), the
name of a projection acquires (if the move is done from right to left)
or loses (if it is done from left to right) a pair $(i,i \pm 1)$. In addition,
all other numbers that are larger or equal to $i$, are going to increase or
decrease by $2$. Since the name to be kept must possess as few pairs as
possible, an appearance of any pair of the form $(i,i \pm 1)$ in a name means
automatic elimination.
\par \medskip Similarly by performing a ``second" Reidemeister move (Figure 4b)
one adds or removes two pairs of the form $(i,j)$ and $(i \pm 1, j \pm 1)$,
while
increasing or decreasing some of the other numbers by $2$ or by $4$.
Therefore the appearance of two pairs $(i,j)$ and $(i \pm 1, j \pm 1)$ in a
name will result to its elimination.
\par \medskip On the other side the ``third" Reidemeister move (Figure 4c) does
not alter the number of pairs; in fact it keeps all but three of the pairs
the same, while replacing pairs $(i,j)$, $(i',k)$ and $(j',k')$, where
$|i'-i|=|j'-j|=|k'-k|=1$, with the pairs $(i,k')$, $(i',j)$ and $(j,k)$.
Therefore the ``minimum pair" criterion is not sufficient and we shall have to
find some additional one.
\par \medskip Let $N$ be a name that does not possess any pair $(i, i\pm 1)$,
neither does it possess any two pairs $(i,j)$ and $(i \pm 1,j \pm 1)$. Let
also $N_1$, $N_2$,...,$N_f$ be the names generated by performing series of
third moves. ($f$ is definitely finite since it cannot exceed
$(2n)! \over n!$). If any of the $f$ names possesses a subset of either the
form $\{ (i,i \pm 1) \}$, or $\{ (i,j), (i \pm 1,j \pm 1) \}$, then not only
this name, but all $N$, $N_1$, $N_2$,...,$N_f$ have to be eliminated, since an
equivalent name with fewer pairs has already appeared at the output. The
``tie breaker" criterion is therefore going to be the following.
\par \medskip Let $N$ be some name and let $|a_i|$ denote the ``companion" of
$i$ with respect to $N$, $a_i$ being positive if $i$ is at the left and
negative if $i$ is at the right; in other words, $j=a_i \Rightarrow (i,j) \in N
\vee (-j,i) \in N$. Let now
$$A(N)=\min_{1\le i\le 2n}|a_i-i|$$
$$B(N)=2n-\max_{1\le i\le 2n}|a_i-i|$$
$$C(N)=\min_{1\le i\le 2n}|a_{i+1}-a_i|$$
$$D(N)=2n-\max_{1\le i\le 2n}|a_{i+1}-a_i|$$
\par \medskip Finally define $M(N)=\rm {min} \bigl(A(N),B(N),C(N),D(N)\bigr)$.
$M(N)=1$ means that $N$ can lose one or more pairs after performing a first or
a second move, while $M(N_r)=1$ for some $r \in \{ 1,2,...,f \}$ implies that
$N_r$ can lose pairs and so can $N$, when third moves are performed to bring
it to $N_r$, which are followed then by a first or a second move. It is thus
desirable to keep
the name whose $M$ is the least and check whether this $M$ is equal to 1. If
so, $N$ is eliminated, otherwise one can be sure that neither this name nor any
of the other equivalent names is reducible.
In case however
where the minimum value of $M$ belongs to more than one names we shall have to
resort to the lexicographical rule.
\par \medskip The process described until now is finite, since given some
name $N$, it only involves names whose order does not exceed that of $N$.
It is possible however that two projections are equivalent, but in order to
connect them one has to perform moves that, initially at least, increase the
number of double points and thus the order of the names involved. Such a
process is more involved than what a first look may suggest, and
the discussion shall be left for section 9 when the relevant algorithm is going
to be
described in detail.
\bigskip
\centerline {6. \qquad TO NAME THE IMPOSSIBLE KNOT}
\bigskip
Choosing just one name for each class and eliminating all others is only part
of the problem faced while trying to devise an algorithm. An additional
difficulty arises from the fact that not all possible ``names" actually
correspond to some knot class. In fact if one takes the first $2n$ natural
numbers, makes $n$ pairs of them and considers the set of these pairs as a
name, one may not be able to draw the resulting projection.
When $n$ goes
in fact to infinity, the probability of getting a drawable name tends to zero.
\par \medskip To understand more clearly this statement, take as an example
the set $\{ (1,3),(2,4) \}$ and try to draw the corresponding projection
(Figure 8). While one can draw the two double points (1,3) and (2,4), the
knot cannot close without creating a third double point, which of course is
forbidden since it would alter the name; the reason for this problem is that
one ``leg" is
inside while the other is outside a loop. While one may obviously see that
closing the loop is impossible, mathematically this is a direct consequence
of the {\bf Jordan Curve Theorem} which states that all closed lines in $R^2$
or $S^2$ which do not intersect themselves divide $R^2$ or $S^2$ into two
disjoint segments [12]. Although this may seem self evident, it does not
occur for
all two-dimensional surfaces. It does not occur for example for $T^2$ or
$RP^2$.
\par \medskip By using the Jordan Curve Theorem one may prove that any name
containing one or more pairs whose elements are either both odd or both even,
does not correspond to any actual knot class. The proof goes as follows. Let
$(i,j) \in N$ as shown in Figure 9. If the number of double points between $i$
and $j$ is odd, it is impossible to draw the name, since the part of the knot
not shown in Figure 9 must enter the loop $i \rightarrow j$ as many times as it
exits and therefore the number of double points between $i$ and $j$ must be
even. In Figure 9 the double points between $i$ and $j$ are numbered
consecutively along the loop, and if they are even then $j-i$ must be odd,
which means that $i$ and $j$ cannot be both odd or both even. It may happen of
course that some numbers of the loop $i \rightarrow j$ are paired together, in
which case
Figure 9 has to be replaced by Figure 10. Such a change does not affect the
argument above, since the number of such points is even, and therefore the
number of the remaining points of the loop must also be even.
\par \medskip Strange as it may appear, this constraint
actually facilitates the devising of an
algorithm. First, the number of possible names for a given number of pairs $n$
is reduced from $(2n)! \over n!$ to $2^n n!$, since $1$ has $n$ even numbers
to choose from as a ``partner", $3$ has $n-1$ choices, $5$ has $n-2$ etc,
while for each pair there are two choices: (odd,even) or (even,odd).
Therefore one may substitute Table II with the following one.
\medskip \centerline {T A B L E \quad III} \medskip
\vbox{\catcode`\*=\active \def*{\hphantom{0}}
\offinterlineskip
\halign{\strut#&\vrule#\quad&\hfil#\hfil&\quad\vrule#\quad&\hfil#\hfil&\quad
\vrule#\quad&\hfil#\hfil&\quad\vrule#\quad&\hfil#\hfil&\quad\vrule#\cr
\noalign{\hrule}
&&\omit \quad $n$ \quad && \qquad $(2n)! \div n!$ \qquad &&
\qquad $2^n n!$ \qquad && \bf Actual Number of Knots &\cr \noalign
{\hrule}
&&0&& **********1&& ********1&& *1&\cr &&1&& **********2&& ********2&&
*0&\cr &&2&& *********12&& ********8&& *0&\cr
&&3&& ********120&& *******48&& *1&\cr
&&4&& *******1680&& ******384&& *1&\cr &&5&& ****30240&& *****3840&&
*2&\cr &&6&& *****665280&& ****46080&& *3&\cr
&&7&& ***17297280&& ***645120&& *7&\cr
&&8&& **518918400&& *10321920&& 21&\cr &&9&& 17643225600&& 185794560&& 49&\cr
\noalign{\hrule}}}
\par \medskip Second, by denoting $f(i)=j$ when $2i-i$ is paired to $2j$ and
$g(i)=0$ or $1$ depending on whether the pair is $(2i-1,j)$ or $(2j,2i-1)$,
ordering the possible names to be studied becomes simpler, since each of them
is actually a pair of an $n$-permutation and an $n$-digit binary number.
\par \medskip The problem of impossible names, however, is still unresolved.
Take for instance the name $\{ (1,4),(3,6),(5,8),(7,10),(9,2) \}$ (Figure 11),
in which every pair contains one odd and one even number. When one tries to
draw
it, things will go smoothly for the first seven points. The eighth point,
however, cannot be paired to the fifth, since one is inside the loop A-B-C-A,
while the other is outside. This is once more a consequence of the Jordan
Curve Theorem. The loops $1-2-3-4$ and $5-6-7-8$ that should intersect at an
even number of points, intersect instead at just one point, namely the $(3,6)$.
We therefore arrive to the following procedure.
\par \medskip Consider any two loops that do not have common segment. (This
condition is necessary, since a common segment will add an infinite number of
common points and thus the notion of an odd or even number of intersection
points will be meaningless). Now consider the points where they intersect,
that is all points $(i,j)$ where $i$ belongs to one loop and $j$ to the other.
If the number of such points is odd for at least one pair of loops, the
corresponding name is impossible.
\par \medskip One should notice of course that not all loops are going to be
of the form $i,i\pm 1, i\pm 2,...,j\mp 1$, where $(i,j) \in N$. Loops may
also be of the form $i,i\pm 1,...,i'\mp 1,i'\equiv j',j'\pm 1,...,j\mp 1$
where $(i,j) \wedge (i',j') \in N$. Such loops possess one or more ``angles".
If two loops share such an ``angle", this should not be counted as an
intersecting point, since at such points the loops ``bounce" instead of
entering or exiting each other (Figure 12).
\par \medskip We shall conclude this section by describing the
lexicographical criterion as promised before. Let $j=f_N(i)$ if $2i-1$ and
$2j$ are paired in some name $N$, and let $g_N(i)=0$ if the odd number $2i-1$
appears on the
left while $g_N(i)=1$ if it appears on the right of its pair. Let now $N_1$
and $N_2$ be two names containing the same number of double points. According
to the lexicographical rule $N_1$ is preferred to $N_2$ if $\exists k \in
\{ 1,2,...,n-1 \} : \forall l<k, \quad  f_{N_1}(l)=f_{N_2}(l) \quad \wedge
\quad f_{N_1}(k)
<f_{N_2}(l)$. In case, however, when $\forall k \in \{ 1,2,...,n-1 \},
f_{N_1}(k)=f_{N_2}(k)$, an additional tie breaker is needed. As such we
choose the following. $N_1$ is preferrable to $N_2$ if $\exists k \in
\{ 1,2,...,n \}:
\forall l<k, \quad g_{N_1}(l)=g_{N_2}(l) \quad \wedge \quad
g_{N_1}(k)<g_{N_2}(k)$. No
further tie breaker is needed, since if for two names $N_1$ and $N_2$,
$f_{N_1}(k)=f_{N_2}(k) \wedge g_{N_1}(k)=g_{N_2}(k) \quad \forall \quad k$,
then $N_1$ and
$N_2$ are identical.
\bigskip
\centerline {7. \qquad PRIME KNOTS - WILD KNOTS}
\bigskip
Let $M_1$ and $M_2$ be two $n$-dimensional manifolds without a boundary. By
removing an $n$-dimensional ball from each of them and by identifying the
resulting boundaries one gets a manifold $M_1 \# M_2$ which is called their
{\bf connected sum}. On the other hand, a manifold that cannot be decomposed
as a connected sum of other manifolds (other than itself and the $n$-sphere
$S^n$) is called a {\bf prime} manifold [13].
\par \medskip Let now $K_1$ and $K_2$ be two knots which are of course
$1$-dimensional manifolds without a boundary. Following the definition above we
may construct their connected $K_1 \# K_2$ sum as shown in Figure 13. Such
knots, which are connected sums of simpler non-trivial knots, do not
usually appear on knot tables, the reason being that once all prime knots are
obtained, one
can easily construct their connected. Let for instance
$N_1$ be the name of $K_1$ having $n_1$ pairs and $N_2$ the name of
$K_2$ with $n_2$ pairs. First one constructs the set $N'_2$ made out of pairs
$(i+2n_1,j+2n_1)$, where $(i,j)$ are the pairs of $N_2$. The name of
$K_1 \# K_2$
is the set $N_1 \cup N'_2$.
\par \medskip This rule can also help us check whether a knot is prime or not.
Let $N$ be the name of a connected sum knot containing whose number of pairs
$n$ is minimum.
Then one can find two numbers $k$ and $l$, at least one of which is
not $1$ or $2n$, which satisfy the following property.
$$\forall (i,j)\in N:k\le i\le l\Leftrightarrow k\le j\le l$$
Therefore through a small modification of the algorithm one can avoid the
appearance of connected sums at the output.
\par \medskip There is another category of knots that shall be briefly
mentioned now,
which are unobtainable; no matter how ``good" is an algorithm and no matter
how high the input parameter, one may never hope to reach the so called
{\bf wild knots}, defined as knots that are inequivalent to all
three-dimensional {\bf polygons} [14].
If one tries to assign a name
to a wild knot following the rules we have established, it is going to contain
an infinite number of pairs, so its turn to appear at the output will never
come. In spite of this fact, the length of a wild knot is finite; the lengths
between two double points of its projection can be, for example, terms of a
decreasing geometric sequence.
\par \medskip Knots which are not wild, and which one can hope to obtain, are
called {\bf tame}.
\bigskip
\centerline {8. \qquad INEQUIVALENT KNOTS}
\bigskip
The discussion until now has mostly been from an ``equivalence" point of view;
in other words, by obtaining all projections with a maximum number of double
points, showing that certain projections are equivalent and by
eliminating the repeated ones, one aimed at keeping just one representative
from each knot class. Such a process, however, can go on indefinitely, since
one may not know {\it a priori} how many moves may be required to show an
equivalence. No matter how many moves one allows for in the program, failing
to prove an equivalence relation does not prove inequivalence.
\par \medskip While knot inequivalence may be proved by using a number of
knot polynomials such as the {\bf Alexander-Conway}, the {\bf Jones} [15],
the {\bf Kauffman} [16] and the
{\bf HOMFLY} [17] ones, such methods face two problems. First, in order to
find a
knot's polynomial one is required to get link polynomials as well, since in
the defining relations $P(a) = \lambda_1 P(b) + \lambda_2 P(c)$ where
$a$, $b$ and $c$ are shown in Figure 14, when
switching a crossing one still gets a knot, but when a double point is
eliminated one gets a link. One would therefore have to generalize the concept
of the ``name" to include links in addition to knots, which would drastically
change and complicate the algorithm. Second, while distinct polynomials mean
inequivalent knots, the converse is not always true; it is possible that
inequivalent knots possess identical polynomials.
\par \medskip Both of these difficulties can be overcome, however, by using
{\bf knot groups} [18]. The group of a knot is defined as the first
homotopic
(fundamental) group of its complement in $S^3$. In order to obtain this group
one considers first all loops in $S^3$ which do not intersect the knot,
regards as equivalent any loops that can be continuously deformed to each
other without intersecting the knot in the process, and thus gets the set of
loop classes, which becomes a group under {\it loop multiplication}. One
interesting point to notice is that while knots are not allowed to have
multiple points in $S^3$, there is no such restriction for the loops of the
knot group.
\par \medskip Two knots that possess distinct groups are of course
inequivalent, but in addition the converse is also (almost) true. Prime knots
with identical groups are equivalent (up to chirality and orientation
reversal) [19]. While in section 7 connected sums were already eliminated,
this
statement provides an additional reason to do so. Comparing knot
groups could be a second method of obtaining knot classes, and were it a
finite process it would be preferrable to the one based on Reidemeister
moves.
\par \medskip First however one has to find the knot groups. The best way
to present a knot group is the {\bf Wirtinger presentation} introduced
by Wirtinger in 1905 [20]. According to this presentation, let $K$ be a knot
and $P(K)$ one of its regular projections on $S^2$
which possesses $n$ double points. The knot group of $K$ is a group
generated by $n$ generators $g_i$, each corresponding to a segment connecting
two undercrossings, satisfying $n$ relations of the form
$g_{i+1}=g_jg_ig_j^{-1}$, where the generators $g_{i+1}$, $g_i$ and $g_j$
start, end and pass through some double point respectively, as shown in Figure
15. While the presentation depends on the projection chosen, the group
itself does not.
\par \medskip Finding presentations of two knot groups does not tell us yet
whether the knots are equivalent or not, unless of course the presentations
are identical. Otherwise, it is possible that two knot groups are the same,
even though their presentations may look quite different.
\par \medskip To illustrate the points above, the knot group of the trefoil
shall be compared to the group of the unknot. Let us first begin with the
trefoil. Let $a$, $b$ and $c$ be the generators of the trefoil group
corresponding to the segments A, B and C as shown in Figure 16. The relations
they satisfy are the following.
\smallskip \centerline {i) \quad $b=c^{-1}ac$ \qquad ii) \quad $c=a^{-1}ba$
\qquad iii) $a=b^{-1}cb$}
\smallskip \noindent One may easily check that one of them is redundant, which
is always the case with one of the $n$ relations of a knot group. The other two
yield the relation $aba=bab$, and thus the trefoil group is the {\bf braid
group} $B_3$, which happens also to be the fundamental group of the
configuration space of three identical particles on $R^2$ [21]. This group
is of
course different from the unknot group which can easily be shown to be the
free one-generator group, which is $Z$.
\par \medskip In general, however, a knot group may not come out to be some
group already known having appeared somewhere else. A proof given by
Reidemeister,
although essentially based on the method shown above, has the advantage that it
can be generalized to prove knot inequivalence in more involved cases.
\par \medskip One may show the trefoil's non-triviality by asking the
following question: ``Is it possible to use three distinct colors to paint a
knot projection, provided that color change can only occur at undercrossings,
and at the neighborhood of an undercrossing either just one color is used, or
all three of them?" First one shows that the answer does not change when a
Reidemeister move is performed, which is demonstrated in Figure 17. In this
figure one considers some projection for which the answer to the question is
negative and thus only one color is allowed which is denoted as color X. If a
Reidemeister move is performed, new segments may arise, but still one may not
use any color other than X and thus the answer remains negative [22].
\par \medskip One can clearly see the relation between Reidemeister's proof
and knot groups by identifying the three colors with the permutations
$\pmatrix {1&2&3\cr 2&1&3\cr}$, $\pmatrix {1&2&3\cr 3&2&1\cr}$ and
$\pmatrix {1&2&3\cr 1&3&2\cr}$ of $P_3$. The constraint that either just one or
all
three colors appear at the neighborhood of an undercrossing can be understood
as the fulfillment of the relations of the knot group. Finally the fact that
all three colors can be used to paint the trefoil but only one for the unknot
implies that trefoil group generators may be assigned to all three permutations
mentioned above, which is the full $\lambda _1=2 \quad \lambda _2=1$
conjugacy class of
$P_3$, while this assignment is not possible for the unknot group. The whole
procedure
does not depend on our knowledge of $Z$ and $B_3$ being different from each
other, not even from our ability to recognize these groups.
\par \medskip One can clearly see now how the ``three color" method can be
generalized. The conjugacy class $\lambda _1=2 \quad \lambda _2=1$ of $P_3$ can
be replaced by any other conjugacy class of $P_n$ for any possible value of
$n$. When comparing two knot projections, if at least one conjugacy class
exists
that can be generated by one knot group but not from the other, the groups
are definitely
distinct and thus the projections belong to inequivalent
knots.
\par \medskip This process is infinite too, since no matter for how many
conjugacy classes two knot groups respond identically, one cannot be sure that
the groups are the same. The best one may do, therefore, is to run
simultaneously two programs, one in order to show knot equivalence by
applying Reidemeister moves and one to show inequivalence by checking
various conjugacy classes of permutation groups, and wait until one of the
programs reaches a verdict.
\bigskip
\centerline {9. \qquad A SUGGESTED ALGORITHM}
\bigskip
Based on the theory already mentioned, we are now going to suggest an algorithm
whose purpose is to list all distinct knot classes which having at least one
projection whose number of double points does not exceed some desired value.
This value is the only input data required. This algorithm, which can be run as
a fortran program, is divided into two parts. In the first, our goal is to
find permissible ``shadows", where by ``shadow" one means a ``flat" knot
projection where one does not distinguish between overcrossings and
undercrossings; in other words, the Figures 3a and 3b are considered
identical. Each such shadow will be named not by a set of pairs, but by a set
of two-element sets, since there is no distinction between $(i,j)$ and $(j,i)$.
Once the permissible shadows are separated from the ``forbidden" ones, we
go on with the second part in order to get the actual knot classes.
\par \medskip Using the notation of section 6, where $j=f(i)$ means that
$2i-1$ is paired to $2j$, we start with the shadow $f(i)=i \quad \forall i \in
\{ 1,2,...,n \}$. First one has to check if this shadow is allowed (actually
it is not). If it is, one proceeds to the second part where
different combinations of $g(i)$ are considered, and once this is over one
moves on to the next
shadow. If it is not, one moves directly on to the next shadow, according to a
procedure that shall now shall be described.
\par \medskip Let $f(i)=j$ be a shadow just studied and from
which one wants to move on. Let $f(n)<f(n-1)<...<f(l+1)$ and $f(l+1)>f(l)$ for
some $l$. The shadow that follows, $f'(i)=j$, satisfies the following
relations: a) $i<l
\Rightarrow
f'(i)=f(i)$, b) $f'(l)$ is the smallest among $f(l+1), f(l+2),...,f(n)$
that is larger than $f(l)$, while c) for $i>l$ one has
$f'(l+1)<f'(l+2)<...<f'(n)$.
Following this procedure, which is equivalent to taking the $n!$ shadows in a
lexicographical order, one eventually arrives to the point where $l=0$ or
$f(n)<f(n-1)<...<f(1)$ which means that $f(i)=n+1-i \quad \forall \quad i$.
Once that point is reached and after that shadow has been studied
(it is forbidden too),
the running of the algorithm is over.
\par \medskip Once some shadow $f(i)=j$ is chosen, it is much simpler from
that point on to switch to the notation $j=h(i) \Leftrightarrow$ $i$ and $j$
are paired. In addition, in order to satisfy the fact that our notation is
mod$2n$, we shall also define $h(i \pm 2n) = h(i)$. Permissible shadows are
going to be the ones satisfying the following conditions.
\par \medskip 1) $ \forall \quad i \in \{ 0,1,2,...,2n \}
\quad h(i) \notin \{ i-1,i+1 \} $ to make sure the knot cannot be simplified
through a first Reidemeister move. (The reason this property must also be
satisfied by $i=0$ is to make sure that $\{ 1,2n \}$ does not belong to the
shadow).
\par \smallskip 2) $\not\exists \quad k,l: \quad k<l \wedge (k\ne 1 \vee l\ne
2n) \wedge (k \le i \le l \Rightarrow k \le h(i) \le l)$, in order to avoid a
connected sum.
\par \smallskip 3) If one compares $h(i)$ to $h'(i)=h(i+k)-k$ or $h(k-i)+i$,
modulo $2n$ of course, then $h(i)$ must be lexicographically preferrable or
at least identical to all $4n$ possible values of $h'(i)$ according to the rule
established at the end of section 6. This will ensure that each projection
will only be considered once.
\par \smallskip 4) We finally have to check whether some name belongs to an
actual shadow or whether it violates Jordan's Curve Theorem. This check is in
fact the most time consuming step of the first part, and is done as follows.
\par \medskip We first choose a value of $k$ that runs from $1$ to $2n$. Then,
for a given $k$, we choose numbers $j_1,j_2,...,j_k$ belonging to
$\{ 1,2,...,2n \}$ which are different from each other. These numbers define a
loop $j_1 \equiv h(j_1),h(j_1)+1,...,j_2-1,j_2 \equiv h(j_2), h(j_2) \pm 1,...,
j_3 \mp 1,
j_3 \equiv h(j_3),...,j_1 \pm 1, j_1 \equiv h(j_1)$. One may notice that the
points $j_i$ are the ``angles" of the loop where it changes direction
following $h(j_i)$ instead of $j_i$. Self intersecting loops are of course
excluded.
\par \medskip Once now one is given two such loops, one first has to make sure
that
they do not share any segment.
For those loops that do
satisfy this property one counts the number of intersecting points, that is
the number of points $(i,j)$ where $i$ belongs to one loop and $j$ to the
other. As explained in section 6, common ``angles" they may share should not
be counted. If the number comes out to be even, the test is considered
successful and one continues with the next pair of loops. If it comes out odd,
the impossibility of the shadow has just been proved and now one moves on to
the next
shadow by replacing $f(i)$ with $f'(i)$. Since in fact we do not care about
the actual number of common points but only whether this number is even or
odd, it is simpler to consider some parameter $L$ equal to $1$,
and each time we
reach a common point, to replace $L$ with $-L$. If at the end $L=1$ the test is
successful, if $L=-1$ the shadow is forbidden. Loops sharing one or more
common segments possess an infinite number of points in common, so it would be
meaningless to carry out this test for them; the result for such a pair is
considered successful in advance.
\par \medskip While this process is extremely large due to the big amount of
possible loops and constraints they have to satisfy, it is definitely finite
so that a shadow's permissibility can be checked after a finite amount of
steps. Once this is over, one moves on to the second part of the algorithm in
order to find out which (if any) out of the $2^n$ possible combinations of
overcrossings/undercrossings of a permissible shadow give rise to distinct
knot classes and should thus appear at the output.
\par \medskip Let us assume now that a certain shadow has passed all four
tests described until now and has thus been found to be permissible. One now
considers the projections possessing this shadow, starting from the case
$g_i=0 \quad \forall \quad i$, where $g_i$ was defined in section 6. This
projection is
an {\bf alternating} one since if one travels along the knot, one encounters
alternatively overcrossings and undercrossings; indeed it is the only
alternating projection of the shadow. Once some projection has been studied
and has either been rejected or printed out at the output, one proceeds by
finding the maximum $i$ for which $g_i=0$, changing this $g_i$ from 0 to 1 and,
if $i \ne n$, $\forall \quad j>i$, change $g_j$ from $1$ to $0$. When one is
through with $g_1=0$, $g_i=1 \quad \forall \quad i>1$, instead of proceeding
with $g_1=1$, $g_i=0 \quad \forall \quad i>1$, one leaves the shadow, returns
to the first part of the algorithm and continues with the next shadow. The
reason for ignoring projections where $g_1=1$ is that as
explained in section 5, one only considers names where $1$ is the left element
of its pair, while according to the definition of section 6, $g_1=1$ would
place $1$ at the right. This condition reduces the number of possible
projections by a factor of $2$.
\par \medskip For a given projection now it will be convenient for all $(i,j)$
belonging to the name to define $p(i)=j$ and $p(j)=-i$, and of course
$p(i \pm 2n)=p(i)$. First one checks whether the projection can be
simplified through a second Reidemeister move. If so, one should be able to
find some $i$ for which $|p(i+1)-p(i)|=1$; in such a case the projection
for which this occurs is rejected, and one moves on with the next one.
\par \medskip If no such $i$ can be found,
one proceeds by performing
third Reidemeister moves. First one checks for any three pairs
$(i,j)$, $(i',k)$ and $(j',k')$ where $|i'-i|=|j'-j|=|k'-k|=1$ and performs
third moves as described in section 5. One also performs third moves on the new
projections one gets, and keeps going until no new projections appear. On each
new projection obtained one carries out the test described in section 5, by
comparing
$M(N_{new})$ to $M(N_{original})$. If this
test is inconclusive, one settles this issue with the lexicographical test of
section 6. One should be careful however to put the new projection in a
lexicographically preferred form by replacing each pair $(i,j)$ with the pair
$(k \pm i, k \pm j)$
for the appropriate value of $k$. This may be necessary because while for the
original projection the name used has been checked to be the preferrable one,
once the pairs $(i,j)$, $(i',k)$ and $(j',k')$ have been replaced by $(i,k')$,
$(i',j)$ and $(j,k)$, the new name may not be the preferrable one.
\par \medskip Even when this test is carried out successfully, one is far from
certain that the projection should appear at the outcome, since it might
happen that two projections are equivalent but the series of Reidemeister
moves that connect them begins with moves that increase the number of double
points. Schematically one shows such a case in Figure 18.
\par \medskip The next test to be carried out is to increase the number of
double points by one by performing a first move as shown in Figure 19. One then
performs third moves and compares all new projections to the initial one. One
may see of course that all of the new projections have a larger number of
double points, unless
one can reduce them by performing a first or a second move at one of the new
projections. One must therefore first compare the number of double points of
the initial and any subsequent projections before moving on to any further
criteria. Fortunately as one can see in Figure 19, the number of the first
moves
with which one may begin is limited; adding a pair $(i,i+1)$ is fruitless
unless $|p(i-1)+p(i)|=1$, since it does not lead to further third moves. Once
this move is performed, one is going
to have $|p(i)+p(i+1)|=|\pm [(i+1)-i]|=1$. If one uses this result to
perform an additional first move as shown in Figure 20, the net result is a
second Reidemeister move as one may see in Figure 21. Since this move is
going to be carried out later where it shall be discussed extensively, this
test is not going to take place at this stage. Therefore the process of
increasing the number of double
points by performing first moves is a finite one.
\par \medskip Finally one increases the number of a projection's
double points by performing second moves. While when performing first moves
one might add an $(i, i\pm 1)$ pair for any value of $i$ and proceed
with third moves, one may not just simply add now pairs $(i,j)$ and
$(i \pm 1, j \pm 1)$ for any possible values of $i$ and $j$, the reason being
that by doing so one may arrive to a name that is impossible to be drawn. In
Figure 22 for example, while the segments $i$ and $j$ can be joined to
perform a second move, that is not possible with $k$ and $l$. Since the
ability to perform a second move
depends only on the shadow of the projection, it is
preferrable to find the pairs of ``neighboring" segments such as $i$ and $j$
during the first part of the algorithm. One may also have ready from the
first part the shadows of the projections after the second move is
performed. Once a particular projection is chosen, for each such pair of
segments there are two possibilities, either $i$ goes ``above" $j$ or it
goes ``below".
\par \medskip In order now to find the pairs of neighboring segments $i$ and
$j$ one
has to check which of the loops found before, while testing the name's
possibility, leave the rest of the knot on the same side of $S^2$. In
other words one is looking for loops that do not intersect any other loop.
All segments constituting such a loop may be paired with each other to
perform a secod move.
\par \medskip Unfortunately however the number of second moves one may carry
out is infinite, and there is no way to tell in advance how many of them may
be needed to connect two equivalent projections. It is at exactly this point
that the algorithm becomes infinitely long. If one tries to make it finite by
cutting it at some point, there is no guarantee that some knot class will
not be repeated.
\par \medskip Using the theory discussed up to this point, we constructed a
fortran program that initially considered all possible one-to-one functions
$p(i)$ as
defined above from $\{ 1,2,...,2n \}$ to $\{ -2n,-2n+1,...,-1,1,2,...,2n \}$,
where odd numbers are mapped to even and vice versa. Subsequently it
eliminated all functions corresponding to knots that could be shown to be
either undrawable or equivalent to preferred knots. Were this process finite,
only distinct knots would survive and appear at the output. In order to keep
the process finite, we allowed for at most one second move ``up". The results
we obtained may be summarized as follows.
\medskip \centerline {T A B L E \quad IV} \medskip
\vbox{\catcode`\*=\active \def*{\hphantom{0}}
\offinterlineskip
\halign{\strut#&\vrule#\quad&\hfil#\hfil&\quad\vrule#\quad&\hfil#\hfil&\quad
\vrule#\quad&\hfil#\hfil&\quad\vrule#\quad&\hfil#\hfil&\quad\vrule#\quad&\hfil
#\hfil&\quad\vrule#\cr
\noalign{\hrule}
&&\omit $n$ && \bf Shadows &&
Shadows$\times 2^{n-1}$ && \bf Program's Knots &&
\bf Actual Knots &\cr \noalign
{\hrule}
&&0&& **1&& ****1&& *1&& *1&\cr &&1&& **0&& ****0&& *0&& *0&\cr
&&2&& **0&& ****0&& *0&& *0&\cr
&&3&& **1&& ****4&& *1&& *1&\cr
&&4&& **1&& ****8&& *1&& *1&\cr &&5&& **2&& ***32&& *2&&
*2&\cr &&6&& **3&& ***96&& *3&& *3&\cr
&&7&& *10&& **640&& *7&& *7&\cr
&&8&& *27&& *3456&& 21&& 21&\cr &&9&& 101&& 25856&& 57&& 49&\cr
\noalign{\hrule}}}
\par \medskip In addition to obtaining the same number of knots to the one
appearing in the literature when $n \le 8$, we verified that these knots are
exactly the same. For $n=9$ we verified the fact that 49 of the knots we
obtained are identical to the ones appearing in the literature, while the
other 8 should be identical to some of the previous 49.
\bigskip
\centerline {10. \qquad A SUGGESTED ALGORITHM - THE SEQUEL}
\bigskip Up to now the algorithm described in section 9 is based on starting
with all possible names that have a given number of pairs and where each pair
has an
odd and an even element. Then one eliminates all names that are either
impossible
or can be shown after a certain number of steps to be equivalent to some other
more preferrable name. Finally one prints out all the names that survive. Such
a process is however incomplete, since it fails to show the inequivalence of
the knot projections it prints out. In this section the previous
algorithm shall be extended by using ideas described in section 8, in order to
prove the
(possible) inequivalence between the surviving projections.
\par \medskip Let $m$ be an integer and $\lambda _1, \lambda _2,...,\lambda _k$
a {\bf partition} of $m$ defined as a decreasing sequence of natural numbers
$\lambda _1 \ge \lambda _2 \ge ... \ge \lambda _k$ whose sum is equal to $m$.
As is well known from group theory, there is a one-to-one correspondence
between each such partition and a
conjugacy class of the $m$ - permutation group $P_m$ [23]. One now tries to map
generators of knot groups to permutations belonging to various conjugacy
classes, noticing in fact that due to the defining relations of knot groups,
its generators are conjugate to each other. We demand of course that the
permutations mapped to various generators fulfill the same relations that
the generators do; this constraint is going to make many such mappings
impossible, so that in a number of cases the permutations one gets do not
generate the whole conjugacy class.
\par \medskip Take for instance the conjugacy class of $P_3$ corresponding to
the partition $\lambda _1 = 2 \quad \lambda _2 = 1$, consisting of the
permutations
$\pmatrix {1&2&3\cr 2&1&3\cr}$, $\pmatrix {1&2&3\cr 3&2&1\cr}$ and
$\pmatrix {1&2&3\cr 1&3&2\cr}$. While one may map each such permutation to a
distinct generator of the trefoil group, only one of these three permutations
may be assigned to the (unique) generator of the unknot group $Z$, and this
permutation by itself cannot generate of course the whole conjugacy class. The
method to be therefore applied is the following. First one attempts to generate
conjugacy
classes in this way, and after each such attempt one assigns a ``YES" or a
``NO" answer for each (surviving)
knot projection and each partition used. This attempt is of course finite
since there is a finite amount of permutations to be mapped to a finite number
of generators. Then one compares the answers for
various knot projections. If there is at least one partition which gives a
different response for two knot projections, these projections have been
proven to be inequivalent. Our final goal is to perform enough Reidemeister
moves and then to use enough partitions so that all surviving projections can
be shown to be inequivalent.
\par \medskip The suggested algorithm, which can either be inserted in the
second part of the one described in section 9, or become an additional third
part, goes as follows. One begins with $m=1$ and proceeds with successive
natural
numbers until one reaches some desired maximum value of $m$. For a given $m$
now one
starts with the partition $\lambda _1 = m$ and proceeds until the
partition $\lambda _1 = \lambda _2 = ... = \lambda _m = 1$ is reached, at which
point one
increases the value of $m$ by one. For a given value of $m$ now the process
goes as follows.
\par \medskip Let $\lambda _1, \lambda _2,...,\lambda _k$ be a partition been
studied and which is to be replace with a new one. There is
definitely going to be some number $\sigma$ such that $\lambda _1 \ge \lambda
_2
\ge ... \ge \lambda _{\sigma} \ge 2$, while $\lambda _{\sigma + 1} = \lambda
_{\sigma + 2} = \lambda _k = 1$. One keeps the first $\sigma - 1$ terms
unchanged, replaces $\lambda _{\sigma}$ with $\lambda _{\sigma} -1$ and then
makes $\lambda _{\sigma + 1} = \lambda _{\sigma + 2} = ...$ equal to the new
value of $\lambda _{\sigma}$. One considers as many of these $\lambda $ 's as
possible so that the total sum $\lambda _1 + \lambda _2 + ...$ does not exceed
$m$. If there is some quantity remaining which is less than the new value of
$\lambda _{\sigma}$,
one takes a last $\lambda _f$ equal to this remaining.
\par \medskip Once now a partition has been chosen, all segments starting and
ending at undercrossings are considered, which correspond to the $n$
numbers appearing
at the right of each pair. One maps them to permutations of the conjugacy
class under consideration and checks whether the knot group relations are
satisfied. If so, one stays with this mapping, and once one has found all such
permutations $p_1,p_2,...,p_s$, one performs all
possible combinations $p_i^{-1}p_jp_i$ and $p_ip_jp_i^{-1}$ in order to get
new such permutations. One similarly continues with the new permutations until
one gets all possible permutations. If their set is the full conjugacy class,
one has reached an affirmative answer, if not, the answer is negative.
\par \medskip One difficulty that arises in this process and which is now
going to be discussed, is the form of the knot group relations. By looking at
Figure 15 one may consider some pair $(i,j)$, the generator $g_j$
corresponding to the segment beginning from the undercrossing $j$, the
generator $g_{j'}$, where $j'$ is the largest number appearing at the right of
a pair which does not exceed $j$, and the generator $g_{i'}$ where $i'$ is the
largest number at the right of a pair which does not exceed $i$. These three
generators should satisfy either the relation $g_j=g_{i'}g_{j'}g_{i'}^{-1}$ if
they look as in Figure 23a, or the relation $g_j=g_{i'}^{-1}g_{j'}g_{i'}$ if
they look as in Figure 23b. Which one however?
\par \medskip For a partition where $\lambda _1 =2$ it hardly matters, since
all corresponding permutations are equal to their inverses; this is in fact
the reason that made Reidemeister's proof of the trefoil's non-triviality
look so simple. For $\lambda _1 > 2$ however the distinction is important.
While looking at the drawing of a projection it may look ``obvious" which of
the two is the case, when only a name is available, which is just a set of
pairs of
numbers, the answer is far from being so clear.
\par \medskip One should first notice that for some starting point
one may choose any of the two relations one wishes, and then the relations
for the other points will be with respect to the starting one. So let us
assume that for the pair $(i,j)$ or $(j,i)$, where $i<j$, one knows which one
is the appropriate relation.
For the next pair, which is either of the form $(i+1,j')$ or of the form
$(j',i+1)$, one
considers a parameter $L$ which shall now be defined, and if $L$ comes out to
be equal to $1$, then the relations for the two pairs are the same, if it comes
out to be $-1$, they are different. We shall now show how to calculate $L$.
\par \medskip One begins with $L=-1$. One multiplies it with a factor of $-1$
if $i$ is an overcrossing and $i+1$ an undercrossing or vice versa, while one
leaves it as it is if both are undercrossings or both are overcrossings. One
then multiplies $L$ by a factor of $-1$ for each number between $j+1$ and
$j'-1$ that is
paired to a number between $i+2$ and $j-1$. The value final value of $L$
determines the relation for $i+1$.
\par \medskip An example is illustrated in Figure 24. Here one begins with
$(i,j)$ and ends with $(j+3,i+1)$. Initially $L=-1$, then the presence of the
pair $(i+2,j+2)$ switches $L$ from $-1$ to $+1$, and finally the fact that $i$
is an overcrossing while $i+1$ is an undercrossing switches $L$ to each final
value $-1$.
\bigskip
\centerline {11. \qquad INVERSE AND SYMMETRIC KNOTS}
\bigskip
In this section a problem avoided until now is going to be discussed, which is
how to find whether some knot is ambient isotopic to its inverse and to
its mirror symmetric. Let $f$ be the continuous one-to-one function from $S^1$
to $S^3$ which defines some knot $K$, such that $f(\theta ) = \vec x \quad \in
\quad S^3=R^3 \cup \{ \infty \}, \qquad 0\le \theta \le 2\pi, \quad
f(0)=f(2\pi )$. The {\bf inverse}
of $K$ is defined by the function $f_{inv}(\theta )=f(2\pi -\theta )$, while
its
{\bf mirror image} by the function $f_{m.i.}(\theta )=-f(\theta )$. In general
a knot may or may not be ambient isotopic to its inverse or its mirror
symmetric, depending on whether a series of Reidemeister moves may connect it
to them. Unfortunately knot group theory is not helpful in proving possible
inequivalence, since knot groups of inverse and mirror symmetric knots are
always identical; the best proof to our knowledge of a knot's chirality
involves HOMFLY polynomials, but does not always work; the knot $9_{42}$ for
example is chiral, even though its chirality cannot be proved through its
HOMFLY polynomial.
\par \medskip Let now $K$ be a knot with a name $N(K)$ and $K'$ its inverse
whose name is $N(K')$. These knots are oriented, so when naming them
one has to choose inverse orientations. If one chooses the same starting point
for both of them, then for each $(i,j) \quad \in \quad N(K)$, a pair $(2n+1-i,
2n+1-j)$ will belong to $N(K')$. They are equivalent if there is a series of
Reidemeister moves that connects them; during the process one is allowed to
change the starting points by replacing the pairs $(i,j)$ of some name with
pairs $(k+i,k+j)$ as described in section 5; one is not allowed however to
change orientation by replacing them with $(k-i,k-j)$.
\par \medskip Unfortunately the chirality problem is not that simple. As
already showed in Figures 7a and 7b, by merely inverting the order of the
pairs of a name and thus substituting $(i,j)$ with $(j,i)$ one does not
necessarily get the mirror image of a knot; similarly, once given the name of
a knot one may draw either the knot or its mirror image and still be
consistent with the name. The reason is that while replacing $(i,j)$ with
$(j,i)$ implies a reflection with respect to the $xy$ plane, a reflection
with respect to a plane passing through the $z$ axis does not alter a name; in
order to get rid of this ambiguity, the first double point
of a knot shall be constrained to look like in Figure 23a. Once this has been
agreed upon, there is
no further ambiguity, unless the knot is a connected sum, and this is one
further reason to consider only prime knots. One is now able to move on to
compare the name $(i,j)$ with the name $(j,i)$.
\par \medskip While the convention above takes care of the ambiguity, one has
to be careful while performing Reidemeister moves; if either the first double
point is involved, or if one has to switch to another first double point in
order to reach a ``preferrable" name, the new first double point may look
instead like the one in Figure 23b, so one may have to change all new pairs
$(i,j)$ with $(j,i)$ in order to avoid performing an undesirable reflection.
Fortunately the rules established in section 10 can help us here too, by
telling us whether a switch from Figure 23a to 23b has occured when the
first double point has changed. As far as Reidemeister moves are concerened,
one had better avoid involving the first double point in such moves, if
necessary by switching the first double point before and after the move.
\bigskip
\centerline {12. \qquad KNOTS ON A LATTICE}
\bigskip
Until now the problem of knot classification was studied by using knot
projections on $S^2$. Knots however are three-dimensional objects, and it
would be interesting if one could study them on the full three-dimensional
space. If one uses as three-dimensional space $S^3$, one faces the continuity
problem; the space of knots is continuous while algorithms always involve
discrete quantities and thus would not be useful in such a case. In this
section we shall elaborate on an idea first suggested by Pippenger,
Janse van Rensburg, Soteros, Sumners and Whittington [24], which is to
approximate
$S^3$ with a cubic lattice $Z^3$ and study knots on this lattice.
\par \medskip  One begins with $Z^3 = \{ (m,n,l) / \quad m,n,l \quad \in \quad
Z\}$
which is the cubic lattice. A tame knot may always be continuously deformed
to a {\bf polygon} whose vertices belong to the lattice, while successive
vertices are neighboring points on the lattice. The defining
function of such a knot is determined by the values of $f({2\pi k\over n})$
for some $n \in N$ and $k \in \{ 0,1,2,...,n-1 \}$. In addition, the following
relations
must be satisfied.
\par \medskip i) \qquad $k\ne k'\Rightarrow f({2\pi k\over n})\ne
f({2\pi k'\over n})$
\par \smallskip ii) \qquad $|\vec {v_k}| \equiv |f({2\pi (k+1) \over n}) -
f({2\pi k \over n})| = 1$
\par \smallskip iii) \qquad
$f({2\pi (k+\epsilon )\over n})=(1-\epsilon )f({2\pi k\over n})+\epsilon
f({2\pi (k+1)\over n})\quad \forall \epsilon \in (0,1)$
\par \medskip
Therefore one may write such a knot as an ordered $n$-ad $a_1,a_2,...,a_n$,
where $a_i \in \{ 1,2,3,4,5,6 \} \forall i \in \{ 1,2,...,n \}$, the
value of $a_k$ determining the direction of $\vec {v_k}$. In particular, the
values $1$, $2$ and $3$ indicate that $\vec {v_k}$ is the unit vector along
the positive $x$, $y$ and $z$ axes, while values of $4$, $5$ and $6$
indicate that it is along the negative axes $-z$, $-y$ and $-x$.
\par \medskip Not all such $n$-ads are suitable of course; there must be as
many $4$'s, $5$'s and $6$'s as $3$'s, $2$'s and $1$'s in order for the knot to
close, that
is to have $f(2\pi )=f(0)$. This statement requires that $n$ is even. On the
other side no $i$ and $j$ other than $0$ and $n$ should share this property,
or in other words one is not allowed to have a subsequence $a_{i+1},a_{i+2},
...,a_j$ with as many $4$'s, $5$'s and $6$'s as $3$'s, $2$'s and $1$'s; if
such a subsequence did exist, one would have $f({2\pi i\over n})=f({2\pi
j\over n})$ and $f$ would not be one-to-one.
\par \medskip Once these two constraints are satisfied however any ordered
$n$-ad will do, and such an $n$-ad will serve as the name of a knot. One
should contrast the simple form of these constraints with the rather
elaborate form of the constraints of the previous algorithm. What one needs
now in order to proceed, is first a new definition of knot complexity to
replace the number of double points, second some name changes (moves) that
will play a similar role to the Reidemeister moves, and third some method
that can show possible knot inequivalence.
\par \medskip One may easily define the ``length" of the knot, $n$, as a
measure of knot complexity, and thus provide such a number at the program's
input; what is expected as output is a set of knot names
$a_1,a_2,...,a_m$ for $m \le n$ such that no two knots corresponding to these
names are ambient isotopic, and all names missing from the output should
be equivalent to names that do. It would be
interesting to see whether the two notions of complexity, namely the number of
double points and the length as defined above, agree with each other; in other
words, if a knot that is ``simpler" to another knot with respect to one
definition is also simpler with respect to the other. For the first
few cases checked this is indeed the case, and the trefoil remains the
simplest non-trivial knot followed by {\it Listing's} knot.
This is very unlikely however
to be the case for all knots.
\par \medskip As far as
equivalence moves are concerned, one gets two such possibilities which are
shown in Figures 25a and 25b; one is not allowed however to perform these
moves if they result to double points, or in other words if one or more of
the points to be occupied have already been occupied before. The name changes
corresponding to such moves are for the first case an ``exchange",
$$a_1,...,a_{k-1},a_k,a_{k+1},a_{k+2},...,a_n \leftrightarrow
a_1,...,a_{k-1},a_{k+1},a_k,a_{k+2},...,a_n $$ and for the second
case a ``pair creation/annihilation"
$$a_1,...,a_{k-1},a_k,a_{k+1},...,a_n \leftrightarrow
a_1,...,a_{k-1},a_j,a_k,7-a_j,a_{k+1},...,a_n$$
\par \medskip (The words ``exchange" and ``pair creation/annihilation" are of
course used metaphorically; let us consider a knot-polygon as a set of $n$
particles and $n$ antiparticles where the value of $u_i$ indicates the
``flavor"
of the $i$-th particle or antiparticle, and let a flavor $-x$, $-y$ or $-z$
denote the antiparticle of the flavor $x$, $y$ and $z$ respectively. Then the
two moves just mentioned may be considered the first as an ``exchange" and the
second as a particle-antiparticle ``creation" (from left to right) or
``annihilation" (from right to left. A knot-polygon may thus be considered as a
loop in the configuration space of a system of three flavor
particles-antiparticles, where particles or antiparticles having the same
flavor are considered indistinguishable. Equivalent knots correspond to
homotopic loops, while the unknot corresponds to the trivial loop. The
implications of this analogy have not yet been studied.)
\par \medskip Among two equivalent names now the one preferred is the
shortest one, and if both are of equal length, the one with the smallest
value of $a_1$; if these too are equal, the one with the smallest $a_2$ and so
on.
\par \medskip By using what we has already been said until now one may write up
a program which initially provides the $6^n$ possible names for some value of
$n$. Then it eliminates names that do not satisfy the necessary constraints.
Finally by using the two
the two moves described above, one checks for possible equivalences and keeps
from each class the preferred name. This process however suffers from the
same deficiancy with the previous algorithm, namely that it is not finite;
when two names are given, one cannot know in advance how many moves may be
required at most to show their equivalence, so failure to prove equivalence
after a certain number of moves does not prove inequivalence. In addition,
we have yet to come up with some procedure such as the one based on knot
groups that directly proves inequivalence. Take for instance the $24$-ad
$a_1=a_2=a_3=1,a_4=a_5=2,a_6=3,a_7=a_8=6,a_9=a_{10}=a_{11}=5,a_{12}=a_{13}=4,
a_{14}=1,a_{15}=a_{16}=2,a_{17}=a_{18}=a_{19}=3,a_{20}=a_{21}=6,a_{22}=5,
a_{23}=a_{24}=4$, which is ``obviously" the trefoil (Figure 26). While this
``looks like" a trefoil and when one projects it in two dimensions one may
even rigorously prove it is a trefoil, just from the number sequence given
above we have not yet obtained a method to show its nontriviality.
\par \medskip Lately a new set of knot characteristics, the {\bf Vassiliev
invariants} has come to our attention, which might be promising in proving
knot inequivalence for three-dimensional algorithms; these invariants are
defined in a larger space than the space of knots, since they include the so
called {\bf singular} knots which are polygons which do intersect themselves.
Extending our concept of name in order to include singular knots is far
simpler than the case of previous polynomials where multicomponent links
had to be included, since now one merely has to disregard the constraint
$i \ne j \Rightarrow f({2\pi i \over n}) \ne f({2\pi j \over n})$. The
dependence however of the invariants on {\bf actuality
tables} makes the problem more involved than what it looks at first sight.
\par \medskip When comparing the two algorithms one arrives to the following
conclusions.
First, the lattice algorithm is much easier to be written down and there are
fewer ``traps" that may cause errors; the projection algorithm is however
much faster to run, since the values of $n$ involved are much smaller. The
trefoil for example possesses only $3$ double points but has a length of $24$.
In addition, the projection algorithm can obtain more reliable results since
it can prove both equivalence and inequivalence.
\bigskip
\centerline {13. \qquad MULTIKNOTS}
\bigskip
In section 2 knots were defined as continuous one-to-one functions from $S^1$
to $S^3$. This definition may be extended to include  $n$-knots or {\bf
multiknots}, which are continuous one-to-one functions from $S^n$ to $S^{n+2}$.
While such knots cannot be visualized in a three-dimensional space, they do
exist as mathematical objects and there has been some interest in their
study [25].
Ambient isotopy may also be defined in a similar manner; Reidemeister moves do
not exist, however, since a possible projection on $S^n$ will possess a
continuously infinite number of double points. As a result, an algorithm based
on such a projection cannot work.
\par \medskip One may do, however, use an extension of the algorithm described
in section 12 based on placing a knot on a (hyper)cubic lattice. We shall now
demonstrate the case for 2-knots which are functions from $S^2$ to $S^4$.
\par \medskip A suitable two-knot for such a purpose consists of unit squares
whose vertices belong to $Z^4$. There are three kinds of equivalence moves
which are shown in Figure 27, and which substitute either one cube
face with the other five (Figure 27a), or two faces sharing a common edge
with the other four (Figure 27b), or three faces sharing a common vertex with
the other three (Figure 27c). (The dots in these figures are the centers of
the substituted faces).
\par \medskip Naming a two-knot however is not so simple as
naming a knot, since the unit squares do not follow one after the other in a
well defined order as was the case with knots. This is because instead
of functions from $S^1$ one has now functions from $S^2$. One solution to this
problem is to replace $S^2$ with the surface of a cube and divide each of the
six faces into $m^2$ squares. The number $m$ is going to be a measure of the
multiknot's complexity.
One may begin with some square which will be
mapped to $(0,0,0,0)$ of $S^4$, and assign a number from $1$ to $8$ to each
pair of
neighboring squares; If one square is mapped to $\vec {r_1}$ and its neighbor
to $\vec {r_2}$, then one requires that $| \vec {r_2} - \vec {r_1} | = 1$,
while the number assigned to their pair will give the direction of
$| \vec {r_2} - \vec {r_1} |$. This assignment of numbers to each pair of
neighboring squares is going to be the name of the knot. Not any name of
course is acceptable, the conditions being now that any loop on the cubic
surface must be mapped to a loop in $S^4$ and vice versa; these conditions
will impose constraints on the name. In Figure 28 for example, the fact that a
loop starting from the center of square A, and via the centers of B,C and D
returns to A, must be mapped to a loop in $S^4$, imposes the constraint
$\{ a_{AB}, a_{BC} \} = \{ a_{AD}, a_{DC} \}$, where $a_{ij}$ is the number
assigned to the neighboring squares $i$ and $j$.
\par \medskip While such an algorithm in principle leaves only the problem
of proving multiknot inequivalence unsolved, in practice it should need an
enormous amount of time to run in order to yield even the first non-trivial
multiknot.
\bigskip
\centerline {14. \qquad CONCLUSION}
\bigskip
The question posed at the title of this paper, that is whether a complete
classification of (tame) knots is possible, continues to defy resolution.
While by using two at least algorithms, one based on the double points of a
two-dimensional projection and one constraining a knot on a cubic lattice, one
may get correctly some of the simplest knots, when running the program for
knots of greater complexity one cannot be sure that the results obtained are
correct and that all knot classes that appear are distinct. Only the first of
these two algorithms in fact provides some help with respect to this matter.
In addition to this theoretical difficulty, the time needed for the algorithms
is enormous so that it is extremely difficult for the moment to improve on the
tables obtained by Thistlethwaite.
\par \medskip It would be interesting to extend the discussion in order to
obtain algorithms that classify multicomponent links; in principle there
shouldn't be any difficulty, although the algorithms would come out much more
involved and the finiteness problem would remain unresolved. A bonus of
obtaining such an algorithm would be its application in calculating knot and
link polynomials.
\par \medskip By using the cubic lattice algorithm an attempt was made to
extend
the discussion to multiknots. In principle this should also be possible, but
the finiteness problem and a number of practical difficulties remain.
\par \medskip We hope that in the future we may return to these questions and
report further progress.
\bigskip
\centerline {A C K N O W L E D G M E N T S}
\bigskip
The author would like to thank Wei Chen for her help at running the fortran
program, locating a number of references and for useful comments about the
manuscript. The author would also like to thank Ms. Klingler for preparing
the figures and professor Weigt for helpful remarks. Finally, the hospitality
and the support of
the Weizmann Institute of Science and of the Brookhaven National Laboratory
during various stages of this work are
extremely appreciated.
\bigskip
\centerline {R E F E R E N C E S}
\bigskip
The references listed here are by no means exhaustive. Inside them the reader
can find further references. For a broader review of the subject of knots the
reader may find especially useful
the books by D. Rolfsen (ref. 6), G. Burde, H.
Zieschang (ref. 9) and L.H. Kauffman (ref.20).
\par \medskip
1) E. Witten, {\it Comm. Math. Phys.} {\bf 121} (1989) 351-399.
\par \smallskip
2) G. Kolata, {\it Science} {\bf 231} (1986) 1506-1508; W.F. Pohl, {\it Math.
Intelligencer} {\bf 3} (1980) 20-27.
\par \smallskip
3) J.W. Alexander, {\it Trans. Amer. Math. Soc.} {\bf 30} (1926) 275-306.
\par \smallskip
4) V. Vassiliev, {\it Translations of Mathematical Monographs} {\bf 28}
(1992); D. Bar-Natan, Ph.D. Thesis, Princeton University, Dept. of Math.,
June 1991; J.S. Birman and X.S. Lin, {\it Inv. Math.} {\bf 111} (1993) 225;
J.S. Birman, {\it Bull. AMS} {\bf 28} (1993) 253.
\par \smallskip
5) K. Reidemeister (1932) {\it Knotentheorie}, Erg. d. Math. 1, No. 1;
H. Seifert, {\it Math. Ann.} {\bf 110} (1934) 571-592.
\par \smallskip
6) K.A. Perko, {\it Proceedings of the A.M.S.} {\bf 45} 262-266;
D. Rolfsen, (1976) {\it Knots and Links}, Berkeley, CA: Publish or Perish,
Inc.
\par \smallskip
7) M.B. Thistlethwaite, (1985) L.M.S. Lecture Notes no 93 pp. 1-76, Cambridge
University Press.
\par \smallskip
8) W. Haken {\it Acta Math.} {\bf 105} (1961) 245-375.
\par \smallskip
9) G. Hemion {\it Acta Math.} {\bf 142} (1979) 123-155.
\par \smallskip
10) K. Reidemeister, {\it Abh. Math. Sem. Univ. Hamburg} {\bf 5} (1927) 7-23.
\par \smallskip
11) G. Burde, H. Zieschang, (1985) {\it Knots}, Berlin: de Gruyter.
\par \smallskip
12) C. Jordan, (1893) {\it Course d' Analyse}, Paris; O. Veblen, {\it Trans.
Amer. Math. Soc.} {\bf 6} (1905) 83.
\par \smallskip
13) J.L. Friedman, R.D. Sorkin, {\it Phys. Rev. Lett.} {\bf 44} (1980) 1100,
{\bf 45} (1980) 148 and {\it General Relativity and Gravitation} {\bf 14}
(1982) 615.
\par \smallskip
14) R.H. Fox, E. Artin {\it Ann. of Math.} {\bf 49} (1948) 979-990; R.H. Fox
{\it Ann. of Math.} {\bf 50} (1949) 264-265; J. Milnor {\it Fund. Math.}
{\bf 54} (1964) 335-338; R. Brode, Diplomatarbeit 1981, Ruhr Universit\"at,
Bochum.
\par \smallskip
15) V.F.R. Jones, {\it Notices of AMS} {\bf 33} (1986) 219-225.
\par \smallskip
16) L.H. Kauffman, {\it Trans. AMS} {\bf 318} (1990) 417-471.
\par \smallskip
17) P. Freyd, D. Yetter, J. Hoste, W.B.R. Lickorish, K.C. Millet, A. Ocneanu
{\it Bull. AMS} {\bf 12} (1985) 239-246.
\par \smallskip
18) L.P. Neuwirth, (1965) {\it Knot Groups}, Princeton University Press.
\par \smallskip
19) B.A. Cipra, {\it Science} {\bf 241} (1988) 1291-1292.
\par \smallskip
20) W. Wirtinger, {\it Jahresber. DMV} {\bf 14} (1905) 517.
\par \smallskip
21) E. Artin, {\it Abh. Math. Sem. Hamburg} {\bf 4} (1926) 47 and {\it Annals
of Math.} {\bf 48} (1946) 101.
\par \smallskip
22) L.H. Kauffman, (1991) {\it Knots and Physics}, World Scientific Publishing
Co. Pte. Ltd.
\par \smallskip
23) A.P. Balachandran, C.G. Trahern, (1985) {\it Lectures on Greoup Theory for
Physicists}, Monographs and Textbooks in Physical Science vol. 3, Bibliopolos,
Napoli.
\par \smallskip
24) N. Pippenger, {\it Discrete Appl. Math.} {\bf 25} (1989) 273-278; E.J.
Janse van Rensburg, {\it J. Phys. A} {\bf 25} (1992) 1031-1042; E.J. Janse van
Rensburg, S.G. Whittington {\it J. Phys. A} {\bf 24} (1991) 5553-5567;
C. E. Soteros, D.W. Sumners, S.G. Whittington {\it Math. Proc. Camb. Phil.
Soc.}
{\bf 111} (1992) 75-91; D.W. Sumners, S.G. Whittington, {\it J. Phys. A} {\bf
21} (1988) 1689-1694.
\par \smallskip
25) P. Cotta-Ramusino, M. Martellini, {\it BF Theories and 2-knots}, preprrint
hep-th 940797.
\end